\documentclass[journal,transmag]{IEEEtran}
\IEEEoverridecommandlockouts
\usepackage{cite}
\usepackage{amsmath,amssymb,amsfonts,bm,upgreek}
\usepackage{tabularx,booktabs}
\usepackage{algorithm,algpseudocode}
\usepackage{graphicx}
\usepackage{textcomp}
\usepackage{xcolor}
\def\BibTeX{{\rm B\kern-.05em{\sc i\kern-.025em b}\kern-.08em
    T\kern-.1667em\lower.7ex\hbox{E}\kern-.125emX}}
\usepackage{fitmc}
\graphicspath{./Figures}

\newcommand{\trans}{\!\top} 

\newcommand{\dndt}{\frac{\mathrm d}{\mathrm dt}}

\DeclareMathOperator{\divergence}{div}
 
\DeclareMathOperator{\curl}{curl} \DeclareMathOperator{\grad}{grad}

\usepackage{draftwatermark}
\SetWatermarkAngle{0}
\SetWatermarkFontSize{12pt}
\SetWatermarkLightness{0}
\SetWatermarkHorCenter{200pt}
\SetWatermarkVerCenter{20pt}
\SetWatermarkText{Copyright~\copyright~2021 IEEE DOI:~10.1109/TMAG.2021.3057828}

\begin{document}
\makeatletter
\newcommand{\linebreakand}{%
  \end{@IEEEauthorhalign}
  \hfill\mbox{}\par
  \mbox{}\hfill\begin{@IEEEauthorhalign}
}
\makeatother

\title{A Two-Step Darwin Model Time Domain Formulation for Quasistatic Electromagnetic Field Calculations\\
\thanks{Corresponding author: M.~Clemens (email: clemens@uni-wuppertal.de)\\This work was supported in parts by the Deutsche Forschungsgemeinschaft (DFG) under grant no. CL143/11-2.}
}

\author{\IEEEauthorblockN{M.~Clemens\IEEEauthorrefmark{1}, F.~Kasolis\IEEEauthorrefmark{1}, M.-L.~Henkel\IEEEauthorrefmark{1}, B.~Kähne\IEEEauthorrefmark{1}, and M.~Günther\IEEEauthorrefmark{2}}
\IEEEauthorblockA{Institute of Mathematical Modelling, Analysis and Computational Mathematics IMACM, University of Wuppertal,}
\IEEEauthorblockA{\IEEEauthorrefmark{1}
Chair of Electromagnetic Theory, Rainer-Gruenter-Stra{\ss}e 21, 42119 Wuppertal, Germany}
\IEEEauthorblockA{\IEEEauthorrefmark{2}
Chair of Numerical Analysis and Applied Mathematics, Gau{\ss}straße 20, 42119 Wuppertal, Germany}}

\maketitle

\begin{abstract}
In the absence of wave propagation, transient electromagnetic fields are governed by a composite scalar/vector potential formulation for the quasistatic Darwin field model. Darwin-type field models are capable of capturing inductive, resistive, and capacitive effects. To avoid possibly non-symmetric and ill-conditioned fully discrete monolithic formulations, here, a Darwin field model is presented which results in a two-step algorithm, where the discrete representations of the electric scalar potential and the 
magnetic vector potential are computed consecutively. Numerical simulations show the validity of the presented approach.
\end{abstract}


\section{Introduction}
\IEEEPARstart{E}{lectromagnetic} field models that do not account for radiation effects are dubbed quasistatic field models. For static fields, the Maxwell equations decouple and enable the consideration of resistive, and capacitive or inductive effects, with either electrostatic, or magnetostatic formulations, separately. For capacitive-resistive effects, the electro-quasistatic (EQS) field model is applicable, while  resistive-inductive effects can be modelled with the magneto-quasistatic (MQS) field approximation \cite{bHausMelcher:01s}.
Quasistatic field  scenarios where inductive, resistive, and capacitive effects need to be considered simultaneously, appear in high-frequency coils and coils of inductive charging systems, where the capacitive effects between the coil windings need to be taken into account and they are a common problem in electromagnetic compatibility, e.g. in automotive engineering.  
For such scenarios, it is quite common to use lumped $R$, $L$, $C$ parameter circuit-type models such as Kirchhoff's model or circuit models in combination with field models used either for parameter extraction or in strong/weakly coupled models, and circuit-formulation oriented partial-element equivalent circuit (PEEC) methods. Rather recently, also field oriented models based on the Darwin field formulation are considered. These quasistatic electromagnetic field models are represented in terms of combined electric scalar and magnetic vector potentials and feature a modified version of Amp\`{e}re's law by eliminating the rotational parts of the displacement currents, i.e., by neglecting the radiation effects in the model.
Darwin formulations \cite{Darwin1920:01s} are not gauge-invariant, and thus, a number of different Darwin model formulations have been considered, \cite{RaviartSonnendruecker1995:01s2}, \cite{Larsson2007:01s}, \cite{KochWeiland2011:01s}, \cite{KochSchneiderWeiland2012:01s}, \cite{inpGarcia2018Chapter1:01s}, \cite{ZhaoTang2019:01s}, \cite{inpBadicsetal2018:01s}, \cite{inpClemensKaehneSchoeps2019:01s}, and \cite{inpKaimori2020:01s}.
The paper is organized as follows. After this introduction, Darwin field models with different established gauge conditions are highlighted. In the third section, a Darwin model is presented which allows to use a two-step numerical solution scheme. Section \ref{sec:numerical} is comprised of numerical experiments with the two-step time domain Darwin formulation, and is followed by conclusions.

\section{The Darwin Field Model}\label{sec:darwin}
Darwin or Darwin-type field models for quasistatic electromagnetic field distributions can be obtained by considering a decomposition of the electric field intensity $\bm{E}$ into an irrotational part $\bm{E}_{\mathrm{irr}}$ and a remaining part $\bm{E}_{\mathrm{rem}}$,
\begin{equation}
\bm{E} = \bm{E}_{\mathrm{irr}} + \bm{E}_{\mathrm{rem}},
\end{equation}
where the irrotational part is represented as the gradient of an electric scalar potential $\varphi$, that is, $\bm{E}_{\mathrm{irr}} = - \grad \varphi$. The remainder part is represented by the time derivative of a magnetic vector potential $\bm{A}$, i.e., $\bm{E}_{\mathrm{rem}} = - \partial\bm{A}/\partial t$. Hence,  
\begin{equation}
\bm{E} = - \frac{\partial}{\partial t} \bm{A} -\grad \varphi\text{,}\qquad
\bm{B} = \curl \bm{A}
\label{eq:E_and_B}
\end{equation}
holds for the electric field intensity and 
for the magnetic flux density, respectively.

The assumption of a quasistatic electromagnetic field model enables the elimination of the rotational parts of the displacement currents,  $\varepsilon \partial^2 \bm{A}/\partial t^2 \cong \bm{0}$, in Amp\`{e}re's law. The result of this elimination is the so-called Darwin-Amp\`{e}re equation 
\begin{equation}
\curl (\nu \curl \bm{A}) + \kappa \frac{\partial}{\partial t} \bm{A} 
                        + \kappa \grad \varphi
                        + \varepsilon  \grad \frac{\partial}{\partial t} \varphi 
                        = \bm{J}_\mathrm{S},
\label{eqn_Darwin-Ampere}
\end{equation}
where $\nu$ is the reluctivity, $\kappa$ is the electric conductivity, $\varepsilon$ is the permittivity, and $\bm{J}_\mathrm{S}$ is a source current density.

The original formulation of the Darwin model \cite{Darwin1920:01s} can be obtained by enforcing a Coulomb gauge $\divergence \left(\bm{A}\right)=0$ corresponding to a Helmholtz decomposition of the electric field intensity $\bm{E}$, i.e., assuming $\curl \bm{E}_{\mathrm{irr}}=0$ and $\divergence \bm{E}_{\mathrm{rem}}= - \divergence \left(\partial\bm{A} / \partial t \right)=0$. Intended to model charges in free space without conductive materials, i.e., $\kappa=0$, $\varepsilon = \varepsilon_0$ and $\nu=\nu_0$, as a consequence to the Helmholtz decomposition, the Gau{\ss} law does not consider the rotational parts of the electric field and yields the electrostatic Poisson equation as a gauge equation. Thus, the original Darwin formulation is given by 
\begin{align}
-\nu_0 \Delta \bm{A} 
                        + \varepsilon  \grad \frac{\partial}{\partial t} \varphi 
                        &= \bm{J}_\mathrm{S},\\
\divergence  \grad \varphi 
                        &= -\rho_{\mathrm{E}} / \varepsilon_0,
\end{align}
which requires to know the electric charge density $\rho_{\mathrm{E}}$ and its motion with $\bm{J}_\mathrm{S} = \rho_\mathrm{E} \bm{v}$ along some velocity vector $\bm{v}$.

To eliminate the free space assumption of the original Darwin model and to include conductors, and permeable and dielectric materials, an application of the divergence operator to the Darwin-Amp\`{e}re equation \eqref{eqn_Darwin-Ampere} results in a modified Darwin field formulation \cite{KochWeiland2011:01s}, \cite{KochSchneiderWeiland2012:01s}, and yields the Darwin continuity equation
\begin{equation}
\divergence \left(\kappa \frac{\partial}{\partial t} \bm{A} 
                        + \kappa \grad \varphi
                        + \varepsilon  \grad \frac{\partial}{\partial t} \varphi\right)
                        = \divergence \bm{J}_\mathrm{S}.
\label{eqn_Darwin-continuity}
\end{equation}
This Darwin continuity equation lacks the radiation term $\divergence (\varepsilon \partial^2 \bm{A}/\partial t^2)$, which is present in the full Maxwell continuity equation. It was shown \cite{KochWeiland2011:01s} that the combined discrete formulation of the 
Darwin-Amp\`{e}re equation \eqref{eqn_Darwin-Ampere} and the Darwin continuity equation \eqref{eqn_Darwin-continuity} results in non-symmetric systems. In addition, the resulting system is singular and requires an additional gauge for the magnetic vector potential in the non-conductive regions of the problem, such as artificial conductivity \cite{KochWeiland2011:01s} or an additional Coulomb-type gauge \cite{ZhaoTang2019:01s}, \cite{inpKaimori2020:01s}
\begin{equation}
\divergence \left(\varepsilon \frac{\partial}{\partial t} \bm{A}\right)=0,
\label{eqn_Coulomb-type_gauge}
\end{equation}
which can be enforced by adding this term with a scaling factor $1/\Delta t$ to the Darwin continuity equation \eqref{eqn_Darwin-continuity},
such that the gauge equation becomes a temporally semi-discrete version of the full Maxwell continuity equation, expressed in terms of the electrodynamic potentials $\bm{A}$ and $\varphi$. The Coulomb-type gauge \eqref{eqn_Coulomb-type_gauge} can be additionally imposed as a third equation via a Lagrange multiplier formulation \cite{ZhaoTang2019:01s}. Both Darwin model field formulations, \cite{ZhaoTang2019:01s} and \cite{inpKaimori2020:01s}, are symmetric and do not require additional regularization.

\section{Two-Step Darwin Model Algorithms}\label{sec:twostep}
The Darwin continuity equation \eqref{eqn_Darwin-continuity} is extended with an additional gauge term $\divergence \left(\kappa \partial\bm{A}/\partial t \right)=0$ \cite{inpClemensKaehneSchoeps2019:01s} to yield the electro-quasistatic equation 
\begin{equation}
\divergence \left(        \kappa \grad \varphi
                        + \varepsilon  \grad \frac{\partial}{\partial t} \varphi\right) = \divergence \bm{J}_\mathrm{S}.
\label{eqn:eqs-continuity}
\end{equation}
The expression $\divergence \left(\kappa\partial\bm{A}/\partial t\right)$ omitted in \eqref{eqn:eqs-continuity} from the Darwin continuity equation \eqref{eqn_Darwin-continuity}, corresponds to explicitly enforcing divergence-free eddy currents in conductive media, i.e., neglecting eventually arising sources and sinks of current densities due to the irrotational parts of the electric field. 
The combination of equation \eqref{eqn_Darwin-Ampere} rewritten as 
\begin{equation}
\curl (\nu \curl \bm{A}) + \kappa \frac{\partial}{\partial t} \bm{A} = 
                        -\kappa \grad \varphi
                        - \varepsilon  \grad \frac{\partial}{\partial t} \varphi 
                        + \bm{J}_\mathrm{S},
\label{eqn_Darwin-Ampere_2}
\end{equation}
with equation \eqref{eqn:eqs-continuity} results in a two-step formulation, where first the electro-quasistatic total current density $\bm{J}_\mathrm{total}= - \kappa \grad \varphi - \varepsilon \grad \left( \partial \varphi/ \partial t\right) + \bm{J}_\mathrm{S}$  
is used as a solenoidal source term with $\divergence \bm{J}_\mathrm{Total} = 0$ to a magneto-quasistatic formulation for the magnetic vector potential $\bm{A}$ represented by the lefthand side of (\ref{eqn_Darwin-Ampere_2}). 
This modified magneto-quasistatic formulation, however, initially does not not address irrotational parts of $\bm{A}$ in the non-conductive regions. While this does not affect the evaluation of $\bm{B}$ in (\ref{eq:E_and_B}), the evaluation of the electric field according to \eqref{eq:E_and_B} involves the expression $\partial\bm{A}/\partial t$ also in the non-conductive regions, which is commonly not covered in magneto-quasistatic field formulations.

To control the irrotational parts of $\bm{A}$, the magneto-quasistatic formulation needs to be regularized. For this, the introduction of a small artificial electrical conductivity $\hat{\kappa}$ in the non-conducting regions has been suggested \cite{KochWeiland2011:01s}. 
In case that $\kappa \gg 1/(\Delta t) \varepsilon$ holds for a given time-step length $\Delta t$, a modified electrical conductivity as e.g. $\hat{\kappa}=\kappa +1/(\Delta t) \varepsilon$ will regularize the formulation, where the resulting time-discrete formulations will feature expressions of the type $1/(\Delta t)\kappa +1/(\Delta t)^2 \varepsilon$ as they occur in second-order time discretization schemes, as e.g. Newmark-beta schemes used for full wave Maxwell-Amp\`{e}re equations \cite{DibbenMetaxas97:01s}.  
Alternatively, a grad-div term augmentation for spatially discretized magneto-quasistatic formulations is applicable  \cite{Bossavit2001:01s}, \cite{ClemensSchoepsDeGersemBartel2011:01s}.

The introduction of a small artificial electrical conductivity $\hat{\kappa}$ in the non-conducting regions is also an established technique to mitigate the static limit instability of the electro-quasistatic formulation (\ref{eqn:eqs-continuity}) that is known to occur for $\frac{\partial}{\partial t} \varphi \rightarrow 0$.

By assuming that $\divergence \left(\kappa\partial\bm{A}/\partial t\right)=0$ holds in equation \eqref{eqn_Darwin-continuity},
the calculation of the electric scalar potential $\varphi$, using \eqref{eqn:eqs-continuity}, is decoupled from that of the magnetic vector potential $\bm{A}$. Thus, it is possible to independently first solve an electro-quasistatic initial-boundary value problem that corresponds to \eqref{eqn:eqs-continuity}, and in a second step solve the modified magneto-quasistatic problem \eqref{eqn_Darwin-Ampere_2}  with the then available total current densities $\bm{J}_\mathrm{Total}$, as depicted in Algorithm \ref{alg1_2-Step_Darwin_TD}.
\begin{algorithm}[ht] 
    \caption{Two-Step Darwin Time Domain} 
    \label{alg1_2-Step_Darwin_TD} 
    \begin{algorithmic}[1] 
        \State Initialize $\varphi(t^0)$ and $\bm{A}(t^0)$;
       \For{$n\gets 0:n_{\mbox{\scriptsize End}}-1$}
            \State Solve problem \eqref{eqn:eqs-continuity} for $\varphi(t^{n+1})$;
       \EndFor
       \For{$n\gets 0:n_{\mbox{\scriptsize End}}-1$}
            \State Solve problem \eqref{eqn_Darwin-Ampere_2} for$\bm{A}(t^{n+1})$;
            \State Evaluate $\bm{E}(t^{n+1})$ and $\bm{B}(t^{n+1})$ with (\ref{eq:E_and_B});
%
       \EndFor
    \end{algorithmic}
\end{algorithm}

Alternatively, it is possible to consecutively execute a solution step for an electro-quasistatic and a magneto-quasistatic field formulation for each discrete timestep, using suitable time stepping schemes \cite{Clemens2005:01s}.
%
%
%
%
%
%

\subsection{Discrete Two-Step Darwin Time Domain Schemes}
Reformulating (\ref{eqn_Darwin-continuity}) and (\ref{eqn_Darwin-Ampere_2}) with a spatial volume discretization scheme, such as the finite
integration technique
\cite{Weiland1996:01s} or the finite
element method with N{\'e}d{\'e}lec elements \cite{Nedelec80:01s}, 
results in a coupled system of time continuous matrix equations 
\begin{equation}
     \Gr^{\trans} \fMkap \Gr               \bm{\upphi} +
     \Gr^{\trans} \fMeps \Gr \dndt         \bm{\upphi} = 
     \Gr^{\trans} \fitvec{j}_{\mathrm{s}}              -
     \Gr^{\trans} \fMkap \dndt             \fitvec{a},
        \label{FIT_Darwin_Continuity}
\end{equation}
\begin{equation}
   \C^{\trans} \fMnu   \C                  \fitvec{a} +
               \fMkap  \dndt               \fitvec{a} =
        \fitvec{j}_{\mathrm{s}} -       \fMkap  \Gr             \bm{\upphi}
   -           \fMeps  \Gr \dndt       \bm{\upphi},
        \label{FIT_Darwin-Ampere_2}
\end{equation}
where $\fitvec{a}$ is the degrees of freedom (dof) vector related to
the magnetic vector potential, $\bm{\upphi}$ is the dof vector of
electric nodal scalar potentials, $\fitvec{j}_{\mathrm{s}}$ is a vector of transient source currents, $\C$ is the discrete curl operator
matrix, $\Gr$ and $\Gr^{\trans}$ are discrete gradient and
(negative) divergence operator matrices. The matrices $\fMnu$,
$\fMkap$, $\fMeps$ are discrete material matrices of possibly nonlinear reluctivities,
conductivities and permittivities, respectively, corresponding to the
specific discretization scheme in use.
Employing e.g. an implicit Euler backward differentiation time stepping scheme with time step $\Delta t$ to \eqref{FIT_Darwin_Continuity} and \eqref{FIT_Darwin-Ampere_2} and $\fM_{\sigma}=\fMkap + (1/\Delta t)\fMeps$ yields a coupled system
\begin{align}
\left[ \Gr^{\trans} \fM_{\sigma} \Gr \right]
                    \bm{\upphi}^{n+1} & =
                    \fitvec{f}_1 (\fitvec{a}^{n+1}),
                    \label{eqn_FIT-BDF1_Darwin_continuity}\\
                \bigl[\C^{\trans} \fMnu \C + \frac{1}{\Delta t} \fMkap \bigr] \fitvec{a}^{n+1} & =
                    \fitvec{f}_2 (\bm{\upphi}^{n+1}),  
                    \label{eqn_FIT-BDF1_Darwin_Ampere}
\end{align}
where, using the notation $\Delta\fitvec{a}^{n+1}=\fitvec{a}^{n+1}-\fitvec{a}^{n}$, the right-hand side vectors are
\begin{equation}
\fitvec{f}_1 = {\Gr}^{\trans} \fitvec{j}_s^{n+1}
                    +
                    \frac{1}{\Delta t}{\Gr}^{\trans} \fMeps {\Gr}\bm{\upphi}^{n}
                    -
                    \frac{1}{\Delta t}{\Gr}^{\trans}\fMkap
                    \Delta\fitvec{a}^{n+1},
                    \label{eqn_rhs_f1}
\end{equation}
\begin{equation}
    \fitvec{f}_2 =
                    \fitvec{j}_s^{n+1}
                    +
                    \frac{1}{\Delta t} \fMkap \fitvec{a}^{n}
                    -
                    \fM_{\sigma} {\Gr} \bm{\upphi}^{n+1}
                    +
                    \frac{1}{\Delta t} \fMeps   {\Gr} \bm{\upphi}^{n}.
                    \label{eqn_rhs_f2}
\end{equation}
System \eqref{eqn_FIT-BDF1_Darwin_continuity}, \eqref{eqn_FIT-BDF1_Darwin_Ampere} are solved for each time step, starting from initial values $\fitvec{a}^0=\fitvec{a}(t^0)$ and $\bm{\upphi}^0=\bm{\upphi}(t^0)$. 
%
Adopting an iterative solution approach with iteration index $i=0,1,2,\ldots$ for each time step $t^{n+1}$ requires to provide an initial guess vector $\fitvec{a}^{n+1}_{i=0}$ with $\fitvec{f}_{1,i=0}=\fitvec{f}_1(\fitvec{a}^{n+1}_{i=0})$.
Inserting the solution vector of \eqref{eqn_FIT-BDF1_Darwin_continuity} rewritten as 
$\bm{\upphi}^{n+1}=\left[ \Gr^{\trans} \fM_{\sigma} \Gr \right]^{-1} \fitvec{f}_{1,i=0}$
into the right-hand side vector equation $\fitvec{f}_{2,i=0} = \fitvec{f}_2(\bm{\upphi}^{n+1}_{i=0})$ of \eqref{eqn_FIT-BDF1_Darwin_Ampere} yields an expression for the next iterative solution vector $\fitvec{a}^{n+1}_{i=1}$. Left application of the discrete divergence operator ${\Gr}^{\trans}$ to this equation using the relation ${\Gr}^{\trans}\C^{\trans}=\bm{0}$ yields the identity ${\Gr}^{\trans} \fMkap \fitvec{a}^{n+1}_{i=1} = {\Gr}^{\trans} \fMkap \fitvec{a}^{n+1}_{i=0}$ and by induction 
\begin{equation}
    {\Gr}^{\trans} \fMkap \fitvec{a}^{n+1}_{i} = {\Gr}^{\trans} \fMkap \fitvec{a}^{n+1}_{0}\quad\forall i\in\{0,1,\ldots\},
\end{equation}
i.e., in exact arithmetics the converged solution $\fitvec{a}^{n+1}$ of the iterative process will maintain the discrete divergence of its initial guess solution $\fitvec{a}^{n+1}_{i=0}$ in the right-hand side \eqref{eqn_rhs_f1} of \eqref{eqn_FIT-BDF1_Darwin_continuity}.
An initial guess $\fitvec{a}^{n+1}_{i=0}= \fitvec{a}^{n}$ yields the result  
\begin{equation}
    {\Gr}^{\trans} \fMkap \fitvec{a}^{n} 
    = {\Gr}^{\trans} \fMkap \fitvec{a}^{0}\quad\forall n\in\{0,1,\ldots\}. 
\end{equation}
Thus, with the choice of the initialization vector $\fitvec{a}^{0}$ of the time integration process at $t^0$ acting as an initial gauge, the difference expression ${\Gr}^{\trans} \fMkap
                    \bigl[\fitvec{a}^{n+1} -\fitvec{a}^{n}\bigr] = \bm{0}$
in \eqref{eqn_rhs_f1} vanishes for all time steps, and thus, system \eqref{eqn_FIT-BDF1_Darwin_continuity}, \eqref{eqn_FIT-BDF1_Darwin_Ampere} gets decoupled.

\section{Numerical Experiments}\label{sec:numerical}
To verify the performance of the proposed two-step Darwin time domain algorithm, two three-dimensional copper coil problems are considered (electrical conductivity $\kappa_{\mathrm{copper}},=5.96\cdot 10^7~\mathrm{S/m}$. Both problems are illustrated in Fig.~\ref{fig:geom}. For each case-study, the computational domain is $\Omega=(\Omega_0\cup\overline{\Omega}_\kappa)\setminus(\Gamma_\mathrm{E}\cup\Gamma_\mathrm{G})\subset\mathbb{R}^3$ and is free from charge and current sources. The bounding surfaces are perfectly conducting, with $\Gamma_\mathrm{G}$ being grounded and $\Gamma_\mathrm{E}$ supplying the transient excitation
\begin{equation}\label{eq:excit}
\varphi(t) = \varphi_\mathrm{max}\cdot f\cdot\min(t,1/f)\cdot\sin(2\pi f t),
\end{equation}
where $\varphi_\mathrm{max}=12~\mathrm{V}$ is the maximum voltage and $f=10~\mathrm{MHz}$ is the excitation frequency with a wavelength $\lambda = 30~\mathrm{m}$ in void. The longest side of the domain $\Omega$ associated with the helical coil is $l=6.3~\mathrm{cm}$, the one associated with the planar coil is $l=1.35~\mathrm{cm}$. Since $\ell\ll\lambda$ in both cases,  the radiation-free assumption of the Darwin field model is justified.
\begin{figure}
\includegraphics[width=0.99\columnwidth]{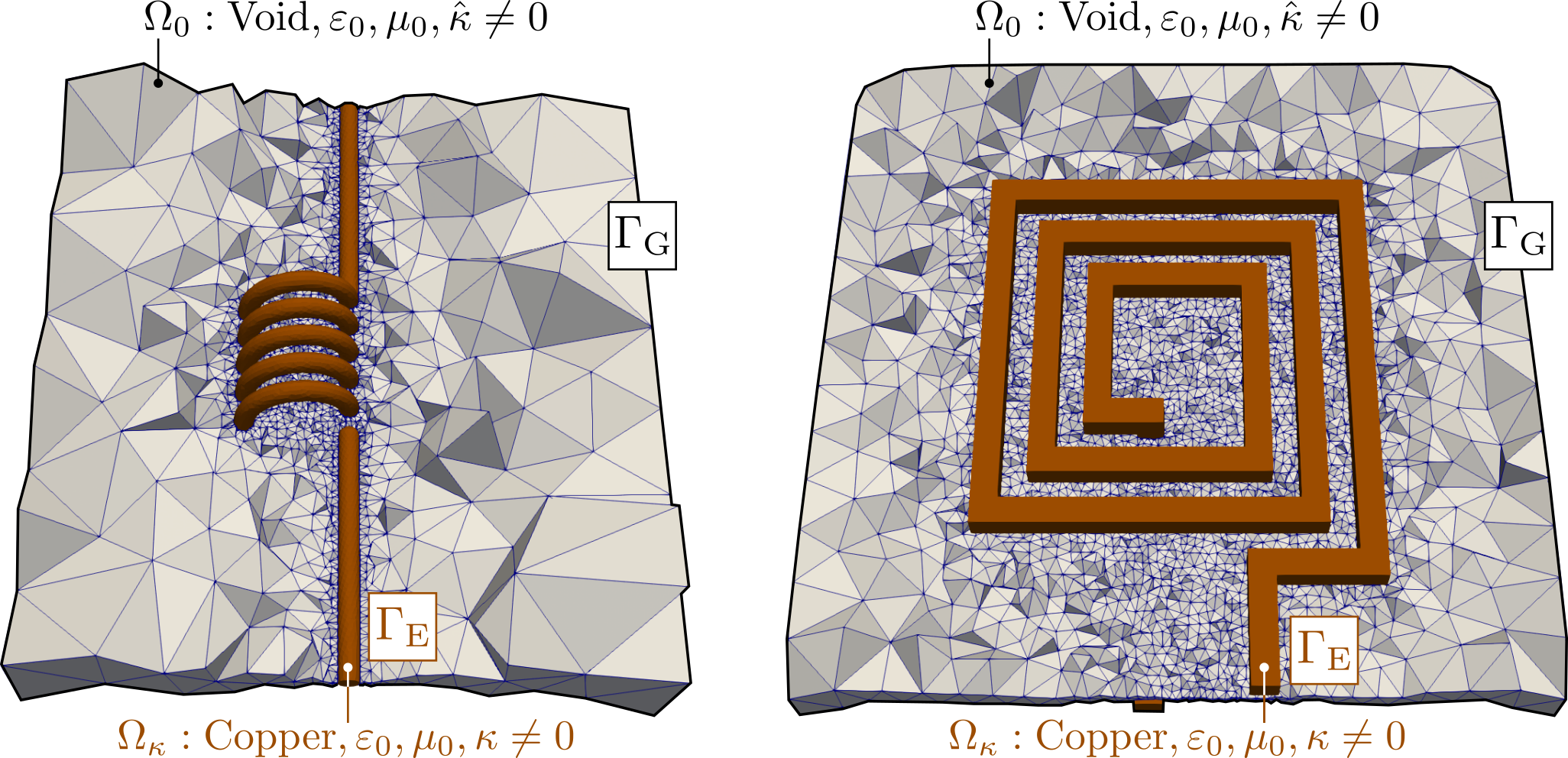}
\caption{In both coil case-studies, $\Omega_0$ is void, with $\hat\kappa$ being a small artificial electrical conductivity ($\hat\kappa=10^{-2}$ S/m), while $\Omega_\kappa$ is occupied by copper.}\label{fig:geom}
\end{figure}

The problems that constitute the two-step algorithm are discretized in space with the FEM, using first-order Lagrange elements for the scalar electric problem and zeroth-order N{\'e}d{\'e}lec elements for the vectorial magnetic problem; see Table \ref{tab:dof} for the number of degrees of freedom in each finite element space. Regarding time-discretization, both problems have been integrated with the trapezoidal method, which is implicit, second-order accurate, and A-stable, with different time steps $\Delta t\in\{2.5,1.25,0.625\}~\mathrm{ns}$ for a total of $t_\mathrm{End}=(n_\mathrm{End}-1)\Delta t=1200~\mathrm{ns}$. With the excitation functions in \eqref{eq:excit}, the two-step algorithm is expected to yield an approximation of a frequency-domain full Maxwell solution, and hence, the latter is used for obtaining reference solutions on the same meshes.
For all linear systems a direct solver is used.
%
\begin{table}\centering
\caption{Number of Degrees of Freedom (dof)}
\label{tab:dof}
\begin{tabular}{lcc}
\toprule
& Lagrange Elements & N{\'e}d{\'e}lec Elements\\
\midrule
Helical Coil & $16\,882$ & $118\,609$ \\
Planar Coil & $41\,357$ & $292\,905$\\
\bottomrule
\end{tabular}
\end{table}

In Fig.~\ref{fig:fields}, the magnetic flux density and the electric field intensity are depicted for each coil. There, the two-step algorithm for the Darwin field model, successfully captures, not only the induction, but also the capacitance between the coil windings. 
Table \ref{tab:freq} depicts the relative differences of the electric field approximations provided by the Darwin field model for the planar coil model at frequencies ranging from $10~\mathrm{kHz}$ up through $1~\mathrm{GHz}$. 
The results show that the remainder part  $-\partial\bm{A} /\partial t$ of the electric field needs to be evaluated in (\ref{eq:E_and_B}) which necessitates the regularization of (\ref{eqn_Darwin-Ampere_2}).

\begin{figure}
\centering
\includegraphics[width=0.8\columnwidth]{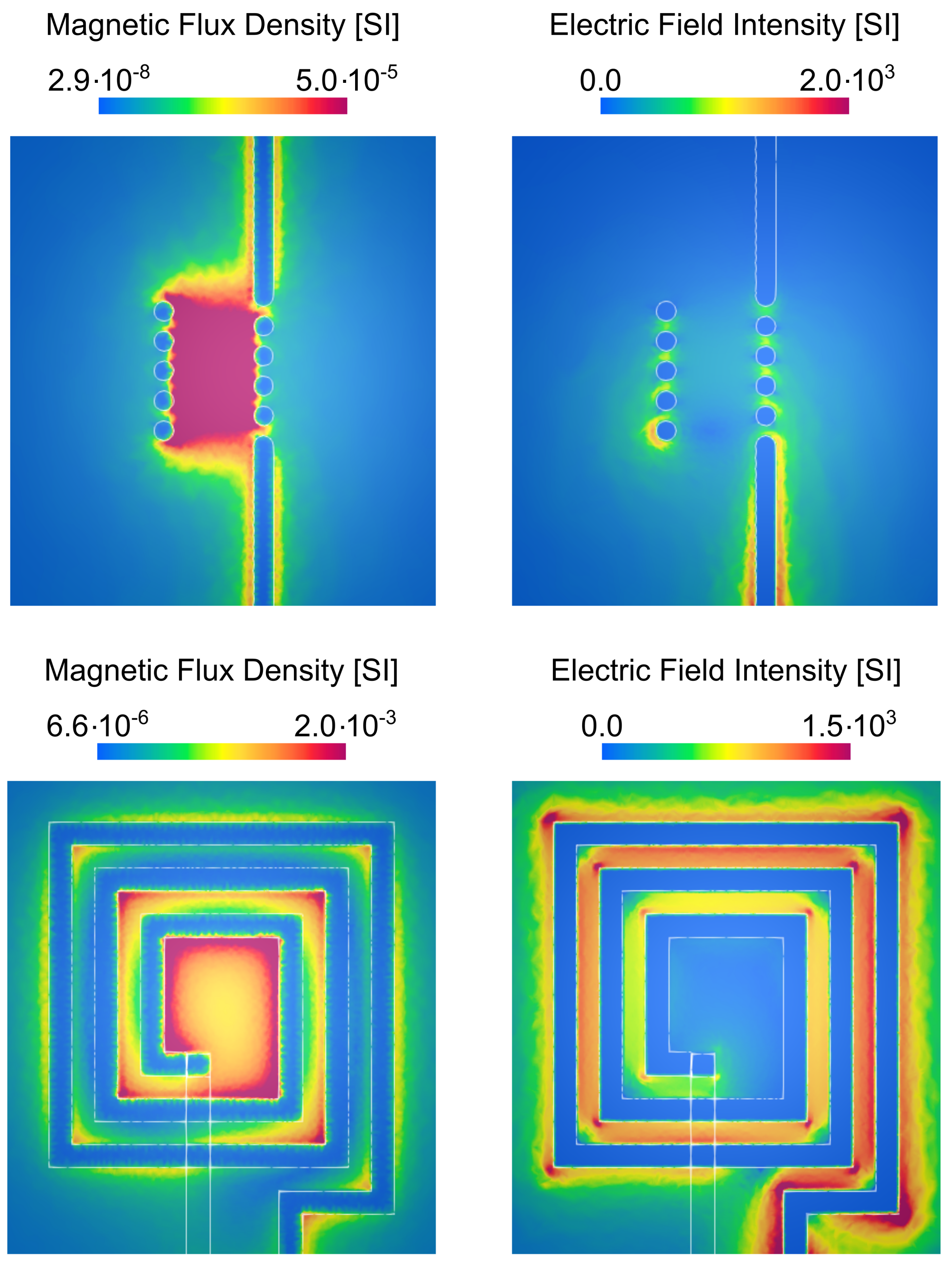}\hfill
\caption{The magnitude of the magnetic flux density and the electric field intensity for the two coil problems at $t=825~\mathrm{ns}$.}\label{fig:fields}
\end{figure}

In Fig.~\ref{fig:error}, the difference between the Maxwell and Darwin field models is quantified with the norm
\begin{equation}
\Vert \operatorname{Re}(\bm{F}_\mathrm{M})-\bm{F}_\mathrm{D}\Vert_{L^2(\Omega)}/\Vert\bm{F}_\mathrm{M}\Vert_{L^2(\Omega)},   
\end{equation}
where $\bm{F}\in\{\bm{B},\bm{E}\}$ is a physical field quantity, computed as in \eqref{eq:E_and_B}, and the subscripts $\mathrm{M}$ and $\mathrm{D}$ stand for the Maxwell and Darwin field models, respectively. In Fig.~\ref{fig:error}, the first row of results is associated with the helical coil, while the second row with the planar coil. In the same figure, the effect of the time discretization scheme is also apparent, with a tendency towards improved accuracy for smaller time steps, since convergence to the time-harmonic solution is expected.
\begin{figure}
\includegraphics[width=0.5\columnwidth]{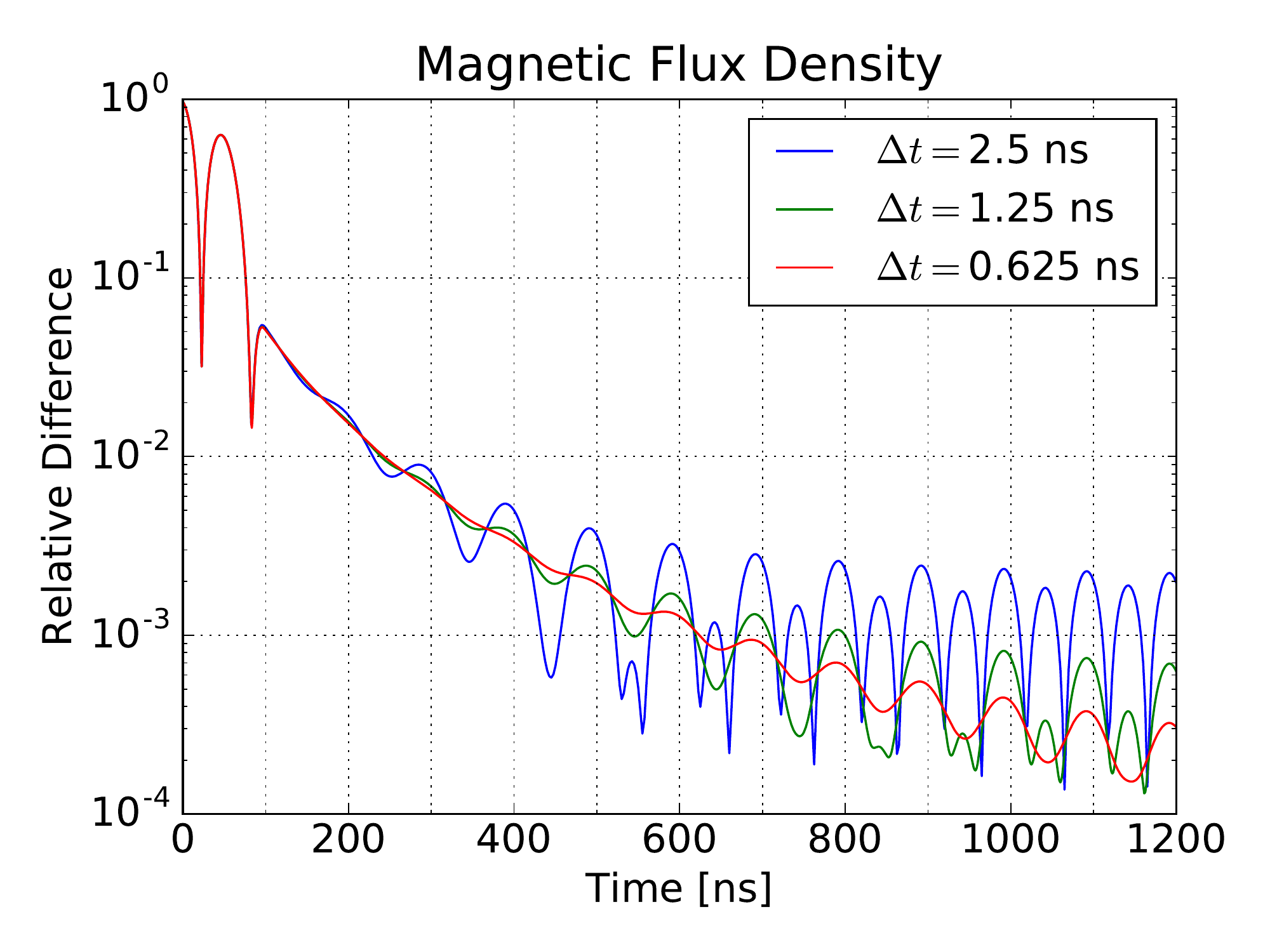}\hfill
\includegraphics[width=0.5\columnwidth]{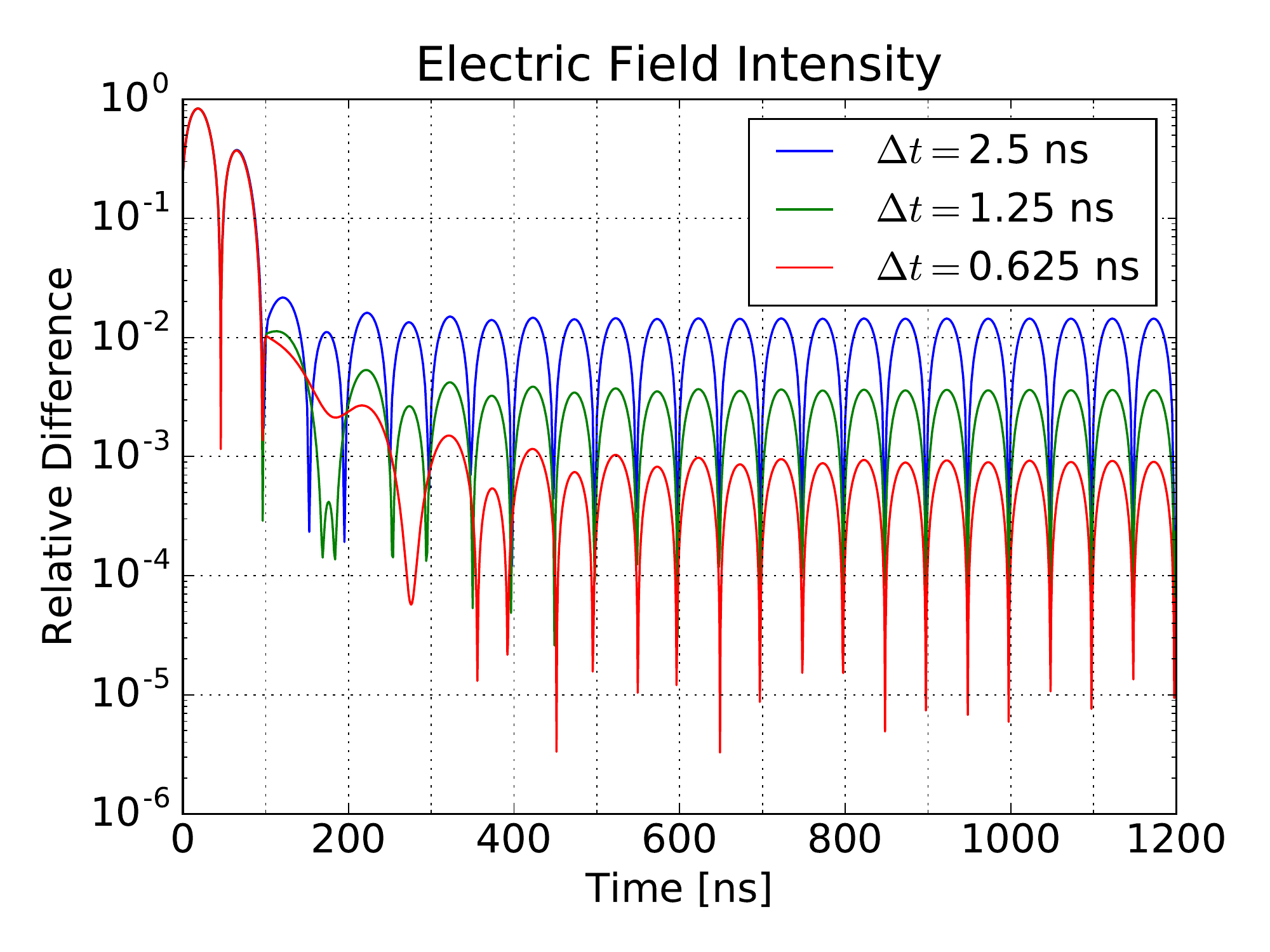}
\includegraphics[width=0.5\columnwidth]{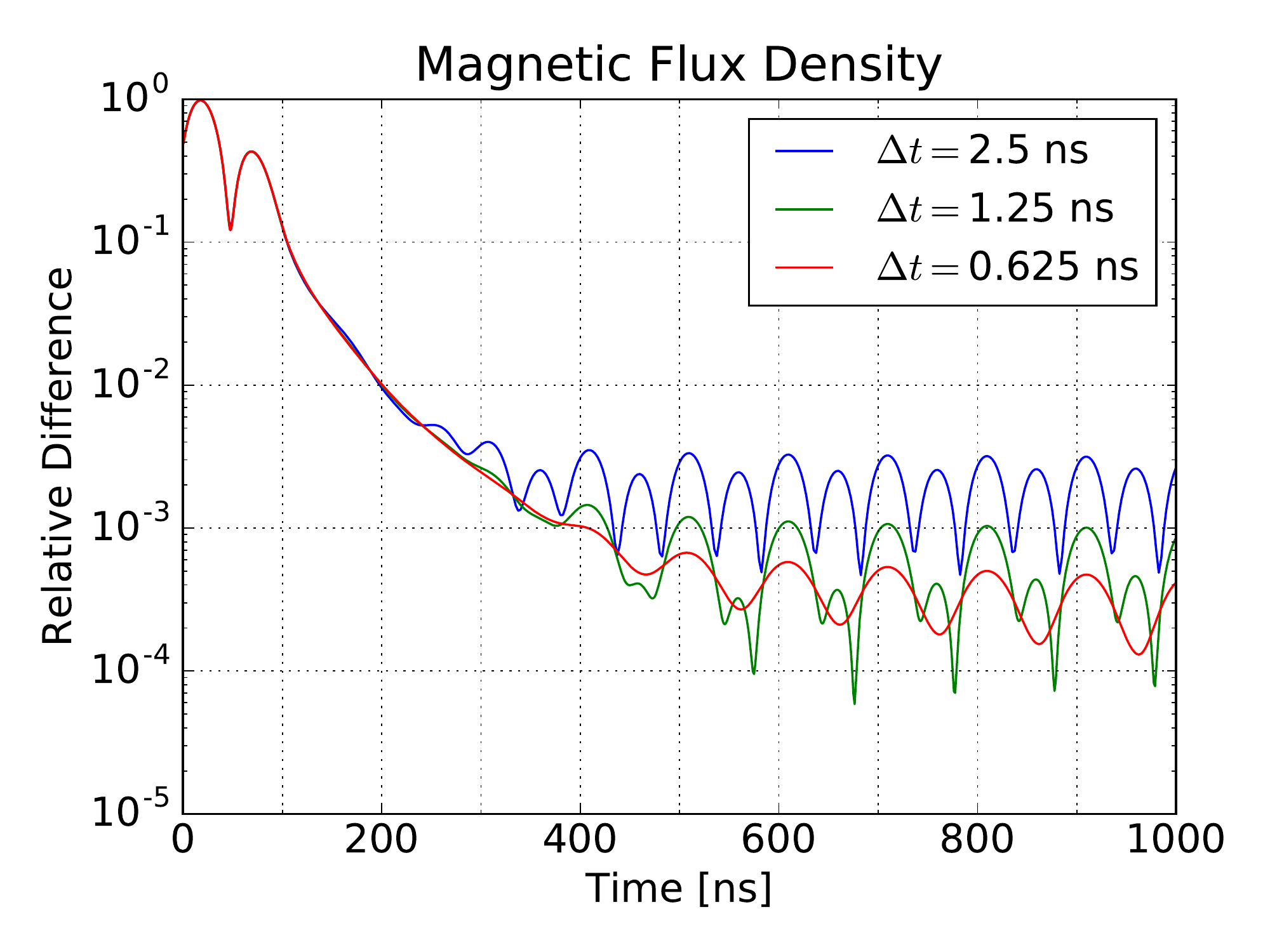}\hfill
\includegraphics[width=0.5\columnwidth]{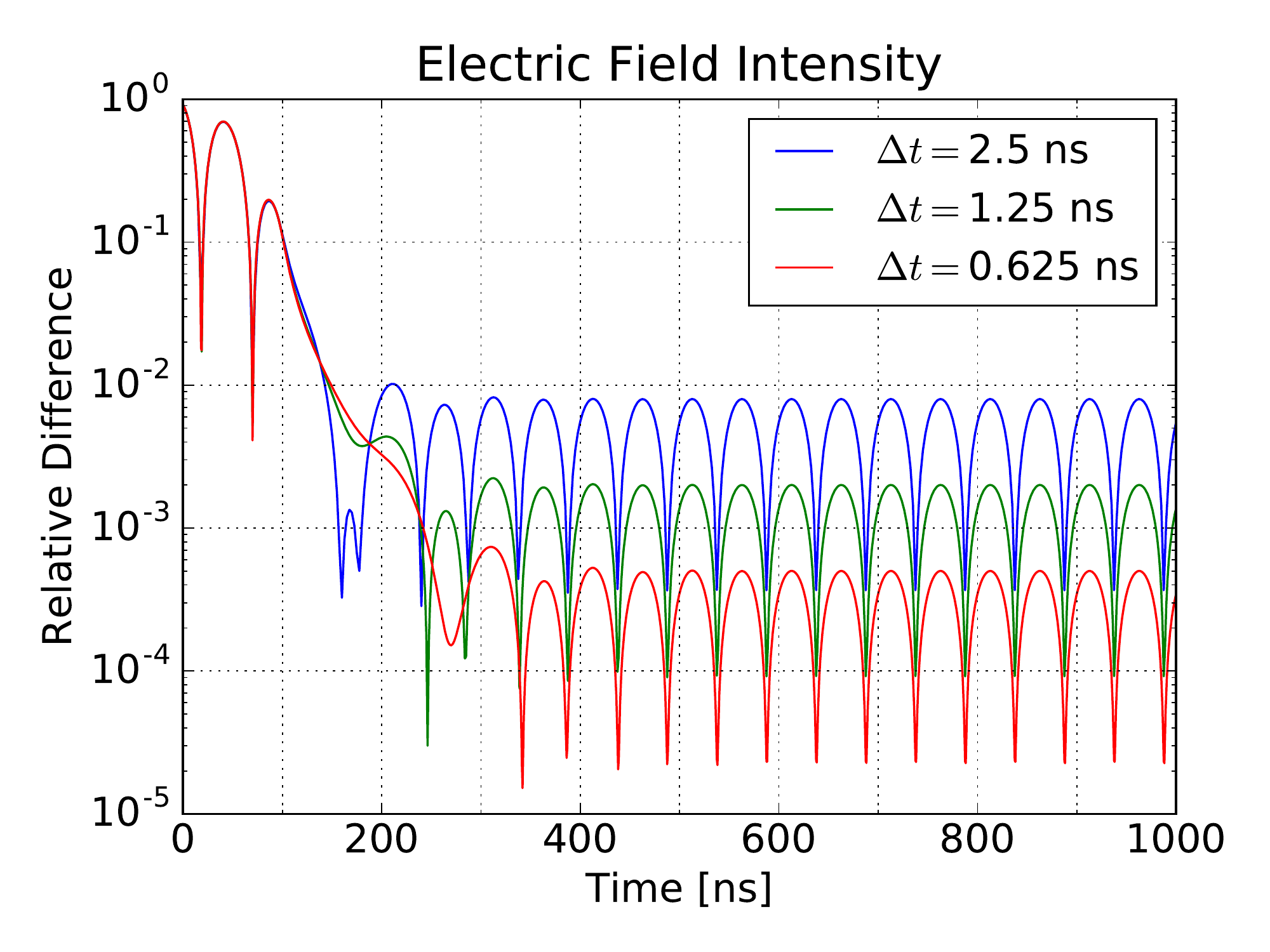}
\caption{The relative difference between the field quantities $\bm{B}$ and $\bm{E}$ for the helical coil (first row) and for the planar coil (second row), computed using a full Maxwell frequency-domain solver and the two-step Darwin algorithm.}\label{fig:error}
\end{figure}
%


\begin{table}\centering
\caption{Relative $\bm{E}$-field Differences at $t=3.125/f$ for the Planar Coil as Functions of the Frequency $f$.}
\label{tab:freq}
\begin{tabular}{ccc}
\toprule
$f$ [Hz] & $\Vert \operatorname{Re}(\bm{E}_\mathrm{M})-\bm{E}_\mathrm{D}\Vert/\Vert\bm{E}_\mathrm{M}\Vert$ & $\Vert \operatorname{Re}(\bm{E}_\mathrm{M})-\bm{E}_\mathrm{D,irr}\Vert/\Vert\bm{E}_\mathrm{M}\Vert$ \\
\midrule
$10^4$ & $8.22\cdot10^{-6}$ & $2.74\cdot10^{-3}$ \\
$10^5$ & $9.31\cdot10^{-6}$ & $2.92\cdot10^{-2}$ \\
$10^6$ & $4.60\cdot10^{-5}$ & $1.59\cdot10^{-1}$ \\
$10^7$ & $7.33\cdot10^{-4}$ & 1.47 \\
$10^8$ & $2.00\cdot10^{-3}$ & 2.28 \\
$10^9$ & $1.41\cdot10^{-2}$ & 2.22 \\
\bottomrule
\end{tabular}

\vspace*{2pt}
\end{table}
\section{Conclusions}
A two-step algorithm for the transient $(\bm{A},\varphi)$ formulation of the quasistatic Darwin field model is introduced and, for the first time, numerically validated against the full system of Maxwell's equations.
The presented two-step Darwin time domain quasistatic field formulation accounts for capacitive, inductive, and resistive effects. The advantages of this scheme result from consecutively combining an  electro-quasistatic and a modified magneto-quasistatic field model,
and thus, it benefits from existing efficient time domain solution techniques that provide flexibility in terms of material non-linearities and excitation profiles.
\bibliographystyle{ieeetr}
\bibliography{refs.bib}
\end{document}